\documentclass[12pt,letter]{article}
\usepackage{amsmath,amssymb,amsthm,amsrefs,graphicx,hyperref}

    \makeatletter
    \newfont{\footsc}{cmcsc10 at 8truept}
    \newfont{\footbf}{cmbx10 at 8truept}
    \newfont{\footrm}{cmr10 at 10truept}
    \renewcommand{\ps@plain}{%
    \renewcommand{\@oddfoot}{\footsc to appear in the electronic journal of combinatorics
  {\footbf 11} (2004), \#R00\hfil\footrm\thepage}}
    \makeatother 

    \title{A Complete Annotated Bibliography of Work Related to Sidon Sequences}

    \author{Kevin O'Bryant\thanks{The author is a National Science Foundation Postdoctoral Research Fellow at the
    University of California at San Diego, grant DMS-0202460.}\\
    \small Department of Mathematics\\[-0.8ex]
    \small University of California, San Diego, USA\\[-0.8ex]
    \small \texttt{kevin@member.ams.org}}

    \date{\small Submitted: May 3, 2004;  Accepted pending revision:  May 17, 2004;
    Published: July sometime, 2004\\
    \small MR Subject Classifications: 11B50, 11B83, 05B10}

\newcommand{\floor}[1]{\left\lfloor #1 \right\rfloor}
\newcommand{\ceiling}[1]{\left\lceil #1 \right\rceil}

\newcommand{\N}{{\mathbb N}}
\newcommand{\Z}{{\mathbb Z}}

\newcommand{\R}{{\mathbb R}}
\newcommand{\C}{{\mathbb C}}
\newcommand{\E}{{\mathbb E}}
\newcommand{\Prob}{{\mathbb P}}
\newcommand{\bigO}[1]{{\cal O}\left(#1\right)}
\newcommand{\littleo}[1]{o\left(#1\right)}
\newcommand{\GF}[1]{{\mathbb F}_{#1}}
\newcommand{\Bhg}[2]{\ensuremath{B^\ast_{#1}{[}#2{]}}}
\newcommand{\Bh}[1]{\ensuremath{B_{#1}}}
\newcommand{\cA}{{\cal A}}
\newcommand{\cB}{{\cal B}}
\newcommand{\cK}{{\cal K}}
\newcommand{\Ruzsa}[3]{{\tt{Ruzsa}}(#1,#2,#3)}
\newcommand{\Bose}[4]{{\tt{Bose}}_{#1}(#2,#3,#4)}
\newcommand{\Singer}[4]{{\tt{Singer}}_{#1}(#2,#3,#4)}

\newcommand{\annotation}[1]{\begin{quotation} \noindent #1 \end{quotation}}
\newcommand{\authorsabstract}[1]{\begin{quotation} \noindent {\bf Author's abstract:} ``#1'' \end{quotation}}
\newcommand{\mathreview}[2]{}%{\begin{quotation} \noindent {\bf Math Review (by #1):} #2 \end{quotation}}
\newcommand{\articlecites}[1]{\begin{quotation} \noindent This article cites #1. \end{quotation}}

\DeclareMathOperator{\ind}{ind}

\newcommand{\MR}[1]{\href{http://www.ams.org/mathscinet-getitem?mr=#1}{{\bf MR~#1}}}

\newtheorem{thm}{Theorem}

\newtheorem{cnj}[thm]{Conjecture}
\newtheorem{dfn}[thm]{Definition}

\begin{document}
\maketitle

\begin{abstract}
A Sidon sequence is a sequence of integers $a_1<a_2<\dots$ with the property that the sums $a_i+a_j$ ($i\le j$)
are distinct. This work contains a survey of Sidon sequences and their generalizations, and an extensive
annotated and hyperlinked bibliography of related work.
\end{abstract}

\section{Introduction}

For a subset $\cA$ of an abelian group (usually $\Z$), define
    $$\cA^\ast(k):=\# \{ (a_1,a_2) \in {\cA} \times {\cA} \colon a_1+a_2=k\}$$
and
    $${\cA}^\circ(k):= \# \{ (a_1,a_2) \in {\cA} \times {\cA} \colon a_1-a_2=k\}.$$

In 1932, Simon Sidon~\cite{1932.Sidon} considered sets of integers with both $\cA^\ast$ and $\cA^\circ$
bounded\footnote{Here is one of his theorems. Suppose that $C$ and $g$ are real numbers. There is a constant
$K=K(C,g)$ such that for every set $\cA$ of positive integers with $\cA^\ast(k)+\cA^\circ(k)\leq g$ (for all
$k>0$) and every sequence $\lambda_a$ with $\sum |\lambda_a|^2\leq C$, there is a function $f:[0,1]\to\C$ which
is bounded by $K$ and $\hat{f}(a)=\lambda_a$ for every $a\in\cA$.}. It is easily shown (and done so in the next
section) that $\cA^\ast(k)\le 2$ for all $k$ if and only if $\cA^\circ(k)\le 1$ for all $k\not=0$. This led
Sidon to ask Erd\H{o}s how large a subset of $\{1,2,\dots,n\}$ can be with the property that $\cA^\ast(k)\leq 2$
(for all $k$)? Since that time such sets, e.g., $\{1,2,5,7\}$, have been known as Sidon sets. In other words,
$\cA$ is a Sidon set if the coefficients of
    $$\left(\sum_{a\in\cA} z^a \right)^2 $$
are bounded by 2.

Unfortunately, many authors call Sidon sets ``$B_2$ sets'', and harmonic analysts use the term ``Sidon set'' to
mean something entirely different. These two factors make searching the literature rather difficult. In Math
Sci-Net, for example, it is impossible to search for a math expression such as ``$B_2$'', and a search for
``Sidon'' returns almost 700 hits, most concerning the harmonic analysts' Sidon sets. Further, there are
several different notations in use, sometimes making it difficult to compare results. For these reasons, I felt
that it would be useful to compile a complete annotated bibliography with a consistent notation and to lay out
the major avenues of research, past, present and possibly future.

In loose terms, requiring a set to have the Sidon property forces it to be thin, e.g., it cannot contain three
consecutive integers. The most basic question is ``How thick can a Sidon set be?'' There are several ways to
make this question explicit, several different settings to explore, and a variety of generalizations. Having
only partially solved these problems, researchers have recently begun turning their attention to the question
``what is the structure of a maximally thick Sidon set?''

In Section~\ref{sec.Terminology}, we define the relevant sets and counting functions and fix the terminology
used in the annotations. In Section~\ref{sec.Constructions}, we present the known general constructions of
(generalized) Sidon sets. In Sections~\ref{sec.Size} and~\ref{sec.Size2} we give the state-of-the-art results
concerning the density of Sidon sets. In Section~\ref{sec.distribution} we state a few of the results
concerning the structure of Sidon sets and their sumsets. In Section~\ref{sec.SpecialSequences} we discuss the
problem of finding a Sidon subset of a given set, such as Sidon sets whose elements are squares or fifth
powers. In Section~\ref{sec.OpenQuestions} we list some of the major unsolved questions.

Finally, we give a partially annotated bibliography of works which either develop or apply the theory of Sidon
sets. The bibliography is, so far as I know, complete and 100\% accurate. The bibliography is heavily
hyperlinked (so you lose something if you print it out), and the links to Math Sci-Net reviews require a Math
Sci-Net subscription. Please email the author regarding any omissions, additions, errors, or clarifications.

\section{Terminology}\label{sec.Terminology}

We begin by defining a generalized Sidon sequence (with parameters $h$ and $g$) to be a sequence $\cA$ such
that the coefficients of
    \begin{equation}\label{gANDh}
    \bigg(\sum_{a\in\cA} z^a\bigg)^h
    \end{equation}
are bounded by $g$. We note that the coefficient of $z^k$ in \eqref{gANDh}, which we denote $\cA^{\ast h}(k)$,
has a number-theoretic interpretation: it is the number of ways to write $k$ as a sum of $h$ (not necessarily
distinct) elements of $\cA$. Also note that this definition is sensible for $\cA$ a subset of any group $G$.

Our notation $\cA^{\ast h}(k)$ is motivated by the notation for Fourier convolution: for any functions $f,g$,
we have
    \begin{equation*}
    f\ast g(k) = \sum_{x\in G} f(x)\overline{g(k-x)}.
    \end{equation*}
Thus, if $\cA$ is the indicator function of the sequence $\cA$ (a common and useful abuse of notation), then
$\cA^{\ast h}(k)$ is exactly the convolution of $h$ copies of $\cA$ evaluated at $k$: $\cA^{\ast h}(k) = \cA
\ast \dots \ast \cA (k)$. For brevity, we write $\cA^{\ast}$ in place of $\cA^{\ast 2}$. Likewise, we hijack
the notation of Fourier correlation:
    \begin{equation*}
    \cA^\circ(k) = \sum_{x} \cA(x) \cA(k+x) = \#\{(a_1,a_2)\in \cA\times\cA \colon a_2-a_1 = k\}.
    \end{equation*}
While convolution is associative, correlation is not. Thus $\cA^{\circ h}$ is ill-defined; we adopt the
convention $\cA^{\circ h}=\cA^{\circ h-1}\circ \cA$.

\begin{dfn}[$\Bhg{h}{g}$ sequence]
A sequence $\cA$ is a $\Bhg{h}{g}$ sequence if the coefficients of $\left( \sum_{a\in \cA} z^a\right)^h$ are
bounded by $g$. If $\cA\subseteq G \not= \Z$, then we call $\cA$ a $\Bhg{h}{g}(G)$ sequence. In particular, if
$G$ is the additive group of integers modulo $n$, then we speak of $\Bhg{h}{g}\pmod{n}$ sequences. We use the
same notation for the property and for the class of sequences with the property, i.e., if $\cA$ is a
$\Bhg{h}{g}$ sequence then we write $\cA\in \Bhg{h}{g}$.
\end{dfn}

Thus, Sidon sequences are exactly the $\Bhg{2}{2}$ sequences. We note that $\cA^\circ$ is bounded by 1 if and
only if $\cA^\ast$ is bounded by 2. For if $\cA^\circ(k)>1$ (with $k>0$), then there are $a_1,a_2,a_3,a_4 \in
\cA$ with $k=a_1-a_2=a_3-a_4$, and at most two of the $a_i$ are equal. This means that
$a_4+a_1=a_1+a_4=a_2+a_3=a_3+a_2$, so that $\cA^\ast(a_1+a_4)\ge 3$ (it is possible that $a_2=a_3$ or
$a_4=a_1$, but not both). Note however that if $\cA=\{2^k,2^k+1\colon k\geq1\}$, then for all $k$,
$\cA^\ast(k)\leq 4$, while $\cA^\circ(1)=\infty$; and if $\cA=\{\pm 2^k \colon k\geq 1\}$ then
$\cA^\circ(k)\leq 3$ for all $k\not=0$, but $\cA^\ast(0)=\infty$. The upshot is that $\Bhg{2}{2}$ sequences are
not only historically important, but they are qualitatively easier to deal with. In essentially all ways, more
is known about Sidon sequences than about $\Bhg{h}{g}$ sequences with $h>2$ or $g>3$.

We use the notation $[n]:=\{1,2,\dots,n\}$. Obviously the $\Bhg{h}{g}$ property is invariant under translation
and dilation, so a supposition of the type ``$\cA\subseteq [n]$'' can usually be replaced with ``$\cA$ is a
subset of an arithmetic progression of length $n$''.

\begin{dfn}[$R$ and $C$]
$R_h(g,n)$ is the largest cardinality of a $\Bhg{h}{g}$ sequence contained in $[n]$. $C_h(g,n)$ is the largest
cardinality of a $\Bhg{h}{g}\pmod{n}$ sequence.
\end{dfn}

\begin{dfn}[$B_h{[}g{]}$ sequence]
A $B_h[g]$ sequence is a $\Bhg{h}{h!g}$ sequence. If $g=1$, then we speak simply of \Bh{h} sequences.
\end{dfn}

We note that many authors define a $B_h[g]$ sequence to be a $\Bhg{h}{h!(g+1)-1}$ sequence (and so a Sidon
sequence, for these authors, is a sequence for which $\cA^{\ast h}(k)\leq 3$). This is not without reason. If
$k=a_1+\dots+a_h$ and the $a_i$ are distinct, then there are $h!$ rearrangements of the $a_i$ which contribute
to $\cA^{\ast h}(k)$. Since there are asymptotically few $h$-tuples from $\cA \times \dots \times \cA$ that
have repeated $a_i$'s, one expects that there is little distinction (asymptotically) between $\Bhg{h}{h!g}$
sets and $\Bhg{h}{h!(g+1)-1}$ sets. This expectation, however, has not been proven to hold except for $h=2$,
$g=1$.

\section{Constructions}\label{sec.Constructions}
\subsection{The greedy algorithm}
The obvious first attempt at constructing a $B^\ast_h[g]$ sequence is to be greedy. Set $\gamma_1=1$, and
define for each $k\ge 1$ the sequence
 $${\cal G}_k=\{\gamma_1, \gamma_2, \dots, \gamma_k\}$$
where $\gamma_{k}$ is the least $m>\gamma_{k-1}$ such that $\{\gamma_1,\gamma_2,\dots,\gamma_{k}\}$ is a
$B^\ast_h[g]$ sequence. Then
 $${\cal G}_h[g]:= \bigcup_{k=1}^\infty {\cal G}_k$$
is an infinite $B^\ast_h[g]$ sequence.

\begin{figure}\label{fig:greedyplot}
\begin{picture}(432,270)
    \put(0,0){\includegraphics{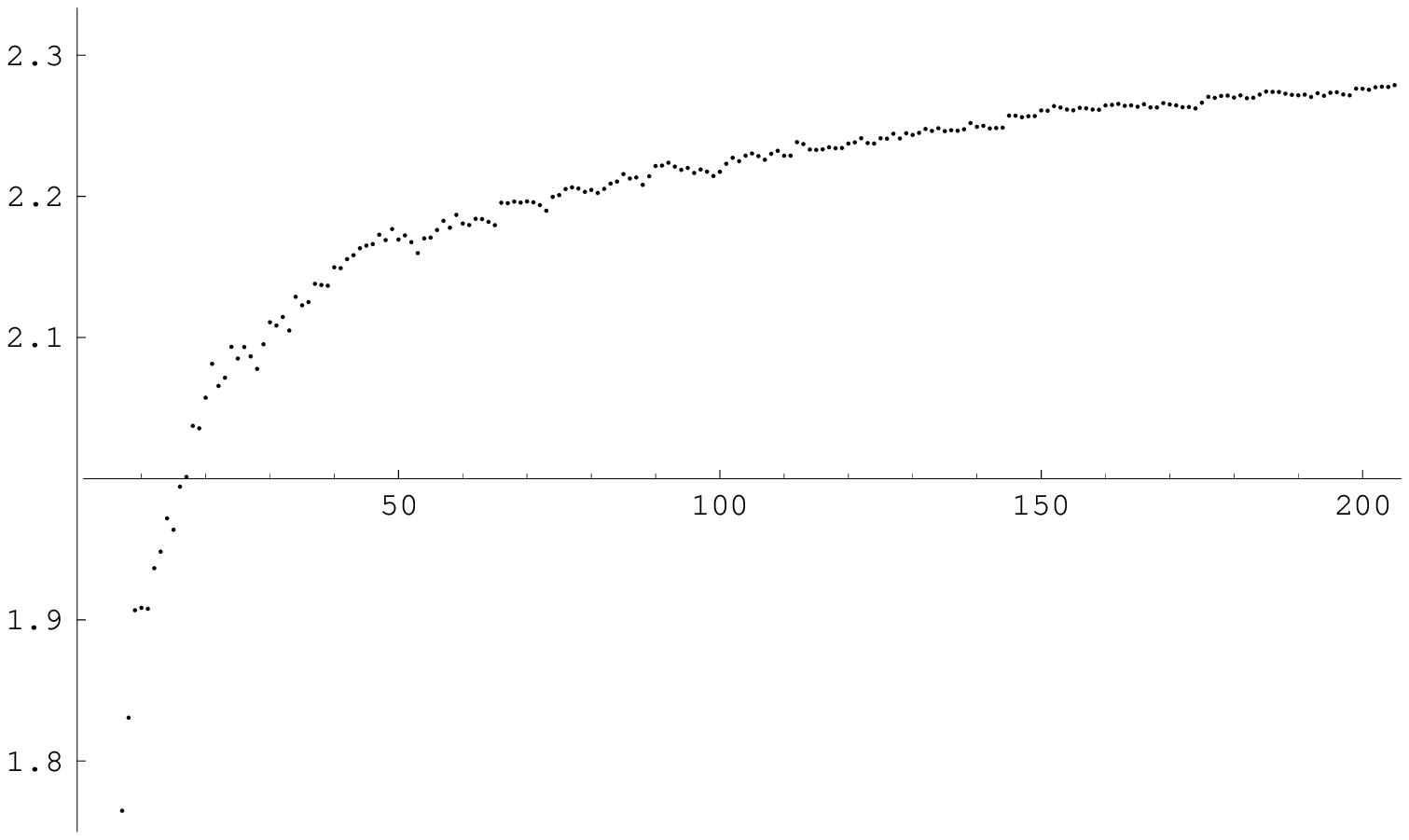}}
\end{picture}
\caption{The points $\left( k, \log_k(\gamma_k)\right)$ for the greedy Sidon set ${\cal G}_2[2]$.}
\end{figure}

Mian \& Chowla~\cite{1944.Chowla.Mian} computed the first terms of ${\cal G}_2[2]$ (i.e., the greedy Sidon
sequence, also called the Mian-Chowla sequence) to be 1, 2, 4, 8, 13, 21, 31, 45, 66, 81, 97, $\dots$. St\"{o}hr
\cite{1955.Stohr} notes that the Mian-Chowla sequence satisfies $\gamma_k < (k-1)^3+1$. Even for this most
simple case, however, the growth of $\gamma_k$ is not well understood. See~\cite{1988.Jia}. The points $\left(
k, \log_k(\gamma_k)\right)$ are shown in Figure~\ref{fig:greedyplot}.

% A straightforward pigeonhole calculation gives
% $$\gamma_k \leq ???.$$

It is interesting to study the greedy sequence with different seeds. Start by setting
$\gamma_1,\dots,\gamma_{r}$, and then continue the sequence in the same manner as above. The resulting sequence
is denoted ${\cal G}_h[g;\gamma_1,\dots,\gamma_r]$.

\begin{cnj}\label{cnj.reciprocals}
For any $h$ and $g$, there are not positive integers $r,\gamma_1,\dots,\gamma_r$ such that
    $$\inf_{\cA \in \Bhg{h}{g}} \left\{ \sum_{a\in\cA} \frac 1a \right\}$$
is achieved with $\cA={\cal G}_h[g;\gamma_1,\dots,\gamma_r]$.
\end{cnj}

This conjecture is dealt with in~\cites{1993.Zhang,2000.Taylor.Yovanof}; it is known that the infimum in
Conjecture~\ref{cnj.reciprocals} (for $g=h=2$) is between $2.16$ and $2.25$.

\subsection{Ruzsa's sets}
A very simple construction of $B_2$ sequences was given by Ruzsa~\cite{1993.Ruzsa}, which is generalized
in~\cite{2004.Martin.OBryant} to $B_2[g^2]$ sequences. Let $\theta$ be a generator of the multiplicative group
modulo the prime $p$. For $k,t\in[p-1]$, let $a_{t,k}$ be the congruence class modulo $p^2-p$ defined by
    \begin{equation*}
    a_{t,k} \equiv t \pmod{p-1} \quad \text{ and } \quad
    a_{t,k} \equiv k \theta^t \pmod{p}.
    \end{equation*}
Define the set
    $$
    \Ruzsa{p}{\theta}{k} := \{a_{t,k} \colon 1\le t < p \} \subseteq \Z/(p^2-p).
    $$
If ${\cal K}$ is any subset of $[p-1]$, then
    $$
    \Ruzsa{p}{\theta}{\cK} := \bigcup_{k\in \cK} \Ruzsa{p}{\theta}{k}
    $$
is a subset of $\Z/(p^2-p)$ with cardinality $|{\cK}|(p-1)$ and
    $$\Ruzsa{p}{\theta}{\cK} \in B_2[|{\cK}|^2].$$

An example is given in Figure~\ref{fig:RuzsaExample}: the row labeled $k$ is $\Ruzsa{13}{2}{k}$. Note that each
row is a translate modulo $13^2-13$ of the row above (and rotated). This is because
    $$a_{t,k\theta}+p\equiv a_{(t+1\bmod{p-1}),k} \pmod{p^2-p},$$
and consequently $\Ruzsa{p}{\theta}{k\theta}+p=\Ruzsa{p}{\theta}{k}$. Since $\theta$ is a generator of the
multiplicative group modulo $p$, this implies that for fixed $p$ and $\theta$ all of the various
$\Ruzsa{p}{\theta}{k}$ are translates of one another.

\begin{figure}\label{fig:RuzsaExample}
    \begin{center}
    $t$\vskip4pt $k$\quad
    \begin{tabular}{|c||cccccccccccc|}\hline
        & 1 & 2 & 3 & 4 & 5 & 6 & 7 & 8 & 9 & 10 & 11 & 12 \\
            \hline \hline
        1&  145& 134& 99& 16& 149& 90& 115& 152& 57& 10& 59& 144  \\
        2&  121& 86& 3& 136& 77&     102& 139& 44& 153& 46& 131& 132 \\
        3&  97& 38& 63& 100& 5& 114& 7& 92& 93&    82& 47& 120 \\
        4&  73& 146& 123& 64& 89& 126& 31&140& 33& 118& 119&    108  \\
        5&  49& 98& 27& 28& 17& 138& 55& 32& 129& 154& 35& 96 \\
        6&  25& 50& 87&    148&101& 150& 79& 80& 69& 34& 107& 84 \\
        7&  1& 2& 147& 112& 29& 6& 103& 128&     9& 70& 23& 72 \\
        8&  133& 110&51& 76& 113& 18& 127& 20& 105& 106& 95&    60 \\
        9&  109& 62& 111& 40& 41& 30& 151& 68& 45& 142& 11& 48 \\
        10&  85&14& 15& 4&     125& 42& 19& 116& 141& 22& 83& 36 \\
        11&  61& 122& 75& 124& 53& 54& 43& 8& 81&     58& 155& 24  \\
        12& 37& 74& 135& 88& 137& 66& 67& 56& 21& 94& 71& 12 \\ \hline
    \end{tabular}
    \caption{Table of $a_{t,k}$ with $p=13$, $\theta=2$. Each row is a Sidon set, the union of any two rows is a
        $\Bhg{2}{8}$ set, the union of any $|\cK|$ rows is a $\Bhg{2}{2|\cK|^2}$ set. Specifically, the row labeled $k$ is
        $\Ruzsa{13}{2}{k}$.}
    \end{center}
\end{figure}

\subsection{Bose's sets}
Bose~\cite{1942.Bose} constructed Sidon sequences by using finite affine geometry. His construction was extended
to $B_h$ sequences (and given in the language of finite fields) in~\cite{1962.Bose.Chowla}, and extended to
$B_2[g^2]$ sequences in~\cite{2004.Martin.OBryant}. Let $q$ be any prime power, $\theta$ a generator of the
multiplicative group of $\GF{q^h}$, and $k \in \GF{q}$, and define the set
    $$
    \Bose{h}{q}{\theta}{k} := \{a \in [q^h-1] \colon \theta^a-k\theta\in \GF{q}\}.
    $$
$\Bose{h}{q}{\theta}{1}$ is $B_h\pmod{q^h-1}$ set, and some work (see~\cite{1945.Bose.Chowla}) has been done on
the question of how quickly these sets can be computed and whether varying $\theta$ (with $h=2$) can help
produce a Sidon set with smaller largest element. If ${\cal K}$ is any subset of $\GF{q}\setminus\{0\}$, then
    $$
    \Bose{2}{q}{\theta}{\cK} := \bigcup_{k\in \cK} \Bose{2}{q}{\theta}{k}
    $$
is a $B_2[|{\cal K}|^2]\pmod{q^2-1}$ sequence.

An example is given in Figure~\ref{fig:BoseExample}, with $q=13$, $h=2$, $\GF{13^2}=\GF{13}[x]/(x^2+2)$, and
$\theta=1+3x$. The column labeled $c_1$ is $\Bose{2}{13}{1+3x\bmod{(13,x^2+2)}}{c_1}$. Note that, as with
Ruzsa's sets, varying $k$ has the effect of translating the set:
    $$\Bose{h}{q}{\theta}{k}=\Bose{h}{q}{\theta}{1}+\log_\theta(k).$$
We note that, unlike Ruzsa's sets, each {\em row} (except the one labeled zero) in Figure~\ref{fig:BoseExample}
is also a Sidon set.

\begin{figure}\label{fig:BoseExample}
    \begin{center}
    $c_1$\vskip4pt $c_0$\quad
    \begin{tabular}{|c||ccccccccccccc|}\hline
         & 0 & 1 & 2 & 3 & 4 & 5 & 6 & 7 & 8 & 9 & 10 & 11 & 12 \\       \hline \hline
        0&  & 77& 147& 21& 49& 35& 91& 7&      119& 133& 105& 63& 161 \\
        1&168& 164& 148& 1& 19& 114& 87& 123&    138& 79&     13& 76& 116 \\
        2&70& 25& 66& 40& 50& 149& 71& 83& 89& 146& 16& 18&     157 \\
        3&112& 23& 58&    108& 125& 31& 92& 20& 67& 113& 60& 82& 131 \\
        4&140&     153& 95& 159& 136& 48& 110& 86& 120& 88& 51& 59&    141 \\
        5&126& 96& 127& 34&     45& 122& 37& 145& 74& 81& 106& 139& 72 \\
        6&14& 90& 93& 137& 128& 15& 10&    130& 27& 152& 101& 33& 162 \\
        7&98& 78& 117& 17& 68& 111& 46& 94& 99& 44&     53& 9& 6 \\
        8&42& 156& 55& 22&    165& 158& 61& 121& 38& 129& 118& 43&     12 \\
        9&56& 57& 143& 135& 4& 36& 2& 26& 132& 52& 75& 11& 69 \\
        10&28&    47& 166&     144& 29& 151& 104& 8& 115& 41& 24& 142& 107 \\
        11&154& 73& 102& 100& 62& 5&     167& 155& 65&    134& 124& 150& 109 \\
        12&84& 32& 160& 97& 163& 54& 39& 3& 30&     103& 85& 64& 80 \\  \hline
    \end{tabular}
    \caption{The least positive integer $k$ such that $(1+3x)^k=c_0 +c_1 x \bmod{(13,x^2+2)}$.
    The column corresponding to $c_1$ is the Sidon set $\Bose{2}{13}{1+3x\bmod{(13,x^2+2)}}{c_1}$.}
    \end{center}
\end{figure}

\subsection{Singer's sets}
Sidon sequences arose incidentally in Singer's work~\cite{1938.Singer} on finite projective geometry. While
Singer's construction gives a slightly thicker Sidon set than Bose's (which is slightly thicker than Ruzsa's),
the construction is more complicated --- even after the simplification of~\cite{1962.Bose.Chowla}. Singer's
construction was extended to $B_2[g^2]$ sequences in~\cite{2004.Martin.OBryant}. No computational work has been
published for Singer's sets.

Let $q$ be any prime power, and let $\theta$ be a generator of the multiplicative group of $\GF{q^{h+1}}$. For
each $\vec{k}=\langle k_1,\dots,k_{h} \rangle \in\GF{q}^{h}$ define the set
    $$
    T(\vec{k} ) := \{0\} \cup \left\{ a \in [q^{h+1}-1] \colon
        \theta^a- \sum_{i=1}^{h} k_i \theta^i \in \GF{q} \right\}.
    $$
Then define
    $$
    \Singer{h}{q}{\theta}{\vec k}
    $$
to be the congruence classes modulo $\frac{q^{h+1}-1}{q-1}$ that intersect $T(\vec k)$. Also define
    $$
    \Singer{h}{q}{\theta}{\cK} := \bigcup_{\vec k\in \cK} \Singer{h}{q}{\theta}{\vec k},
    $$
where $\cK$ is any subset of $\GF{q}^h$. The set $\Singer{h}{q}{\theta}{\langle 1,0,0,\ldots \rangle}$ is a
$B_{h}\pmod{\frac{q^{h+1}-1}{q-1}}$ set. The set $\Singer{2}{q}{\theta}{\langle 1,[k],0 \rangle}$ is a
$\Bhg{2}{2k^2}$ set.

\subsection{Erd\H{o}s \& Tur\`{a}n's sets}
Erd\H{o}s \& Tur\`{a}n~\cite{1941.Erdos.Turan} gave a construction based on quadratic residues. These sets are
substantially thinner than the constructions of Ruzsa, Bose, and Singer given above.

Fix a prime $p$, and let $(k^2)$ be the unique integer in $[p-1]$ congruent to $k^2$ modulo $p$. The set
$\{2pk+(k^2) : 1\leq k <p \}$ is a $B_{2}$ set contained in $[2p+1,2p(p-1)+1]$.

\subsection{Probabilistic Sets}
The seminal paper of Erd\H{o}s \& R\'{e}nyi~\cite{1960.Erdos.Renyi} introducing the probabilistic method to
combinatorial number theory contains a small section on $B_{2}{g}$ sets. The crucial observation made is that
(setting $p_n=1/\sqrt{n}$) the sum
    $$\sum_{n=1}^{N-1} p_n p_{N-n}$$
is bounded independent of $N$. In particular, they prove the existence of an infinite $\Bhg{2}{g}$ set
$a_1<a_2<\cdots$ satisfying
    $$a_k=\bigO{k^{2(1+2/(g-1))}}.$$
See the review of~\cite{1960.Erdos.Renyi} for a more detailed explanation.

See~\cites{1996.3.Kolountzakis, 1998.2.Ruzsa, 1999.Godbole.Janson.Locantore.Rapoport,
2000.Baltz.Schoen.Srivastav, 2001.Ruzsa, Nathanson} for some applications of probability to Sidon sequences, and
vice versa.

\section{The size of finite Sidon sequences}\label{sec.Size}

One of the ``most wanted'' problems is to asymptotically estimate $R_h(g,n)$ for any $h,g$ with $h>2$ or $g>3$.
Also on the ``most wanted'' list is a construction of $\Bhg{h}{g}$ sequences, with $g>h!$, which is dense but is
not merely several $B_h$ sets woven together (such sets do exist:~\cite{1980.Erdos}).

We measure the thickness of $\Bhg{h}{g}$ sets with the quantity:
    $$\sigma_h(g):= \lim_{n\to\infty} \frac{R_h(g,n)}{\sqrt[h]{\floor{g/h!}n}}.$$
Strictly speaking, this limit is not known to exist for $h>2$ or for $g>3$. One should always understand a
lower bound on $\sigma_h(g)$ as being a lower bound on the corresponding $\liminf$, and an upper bound as being
an upper bound on the $\limsup$. If $g=h!$, then we simply write $\sigma_h$.

There are a number of conjectures that are natural to make:
 \begin{enumerate}
    \item   The limit in the definition of $\sigma_h(g)$ is in fact well-defined.
    \item   For each $h$, $\sigma_h(g)$ is an increasing function of $g$.
    \item   For each $h$, $\lim_{g\to\infty} \sigma_h(g)$ is defined and finite.
    \item   If $kh!\le g_1\le g_2 < (k+1)h!$, then $\sigma_h(g_1)=\sigma_h(g_2)$.
 \end{enumerate}
As basic as these questions may seem, they remain unanswered. Regarding question 3, the limit is known to be
bounded between two positive constants. Sadly, there isn't even a conjecture as to growth of $\sigma_h$ as a
function of $h$ (it may well depend on the parity of $h$).

The only explicit values known are $\sigma_2(2)=\sigma_2(3)=1$.

\subsection{$h=2$}
Erd\H{o}s offered USD 500 for answer to the question, ``Is $R_2(2,n)-\sqrt{n}$ unbounded?'' The answer is
likely ``yes''; it is also likely that $R(2,n)-\sqrt{n}$ is nonnegative. Unfortunately, there has been no
progress on these questions since 1941, when it was found that
    $$-n^{\alpha/2} < R(2,n)-\sqrt{n} < n^{1/4}+1,$$
the lower bound holding only for $n$ sufficiently large, with $\alpha$ being a real number such that there is
always a prime between $n-n^{\alpha}$ and $n$ (the current record is $\alpha=0.525$).

\begin{figure}\label{fig:sigma2}
\begin{picture}(432,270)
    \put(0,0){\includegraphics{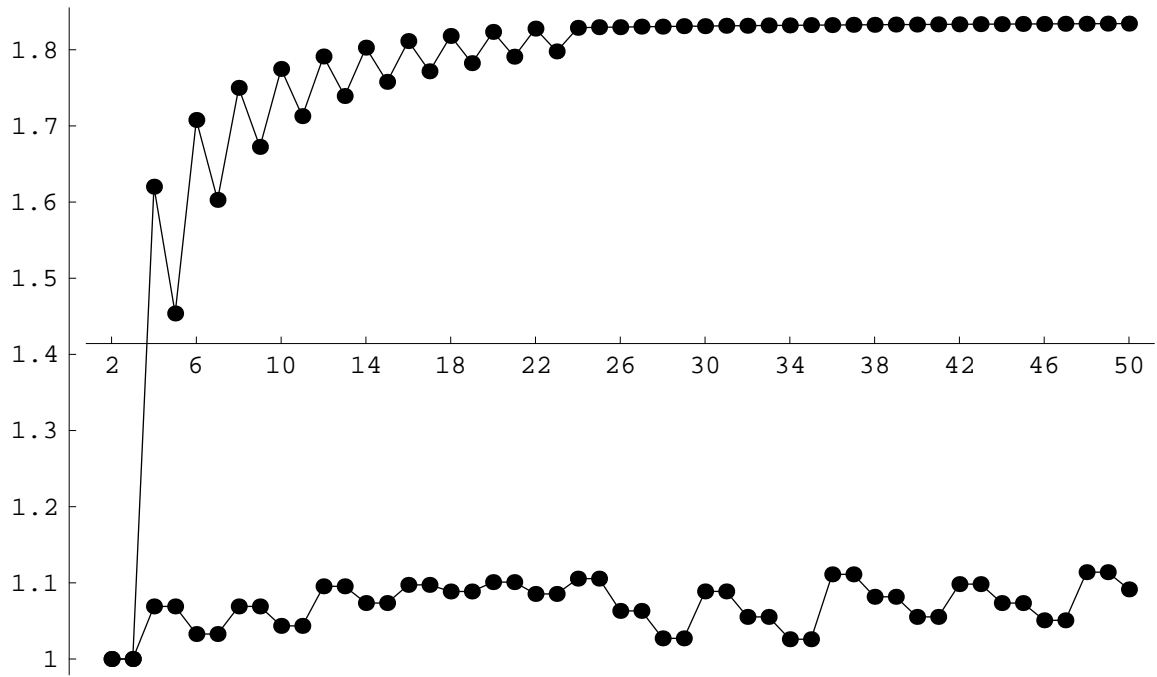}}
\end{picture}
\caption{The best known upper and lower bounds on $\sigma_2(g)$.}
\end{figure}

The best known upper and lower bounds on $\sigma_2(g)$ are shown in Figure~\ref{fig:sigma2}. The lower bound is
$$
\begin{array}{r@{{}\ge{}}l@{{}>{}}l}
    \sigma_2(4)  & \sqrt{8/7}                   & 1.069, \\ \vspace{1mm}
    \sigma_2(6)  & \sqrt{16/15}                 & 1.032, \\ \vspace{1mm}
    \sigma_2(8)  & \sqrt{8/7}                   & 1.069, \\ \vspace{1mm}
    \sigma_2(10) & \sqrt{49/45}                 & 1.043, \\ \vspace{1mm}
    \sigma_2(12) & \sqrt{6/5}                   & 1.095,
\end{array}
\qquad
\begin{array}{r@{{}\ge{}}l@{{}>{}}l}
    \sigma_2(14) & \sqrt{121/105}               & 1.073, \\ \vspace{1mm}
    \sigma_2(16) & \sqrt{289/240}               & 1.097, \\ \vspace{1mm}
    \sigma_2(18) & \sqrt{32/27}                 & 1.088, \\ \vspace{1mm}
    \sigma_2(20) & \sqrt{40/33}                 & 1.100, \\ \vspace{1mm}
    \sigma_2(22) & \sqrt{324/275}               & 1.085,
\end{array}
$$
and for $g\ge 12$
    $$\sigma_2(2g) \ge \sqrt{2} \frac{g+2\floor{g/3}+\floor{g/6}}{\sqrt{6g^2-2g\floor{g/3}+2g}}.$$
In particular,
    $$\lim_{g\to\infty} \sigma_2(g) \ge \sqrt{121/96} > 1.122.$$
The upper bound is a combination of
    $$\sigma_2(2g) \le \sqrt{\tfrac 74 (2-1/g)}$$
and
$$\frac{\floor{g/2}}{g}\sigma_2(g)^2 \leq
\left\{%
\begin{array}{ll}
    1.74043-{1.00483}/g,
    & \hbox{$g\le8$ and even;} \\
    1.58337-\frac{0.026335}g + \sqrt{0.011572-\frac{0.083397}g+\frac{0.00069356}{g^2}},
    & \hbox{$g\ge 10$ and even;}\\
    1.74043-\frac{4.75492}g,
    & \hbox{$g\le 23$ and odd;} \\
    1.58337-\frac{0.071949}g + \sqrt{0.011572-\frac{0.22784}g+\frac{0.0051768}{g^2}},
    & \hbox{$g\ge 25$ and odd.}
\end{array}%
\right.
$$
In particular,
    $$\lim_{g\to\infty} \sigma_2(g) \le 1.839.$$

For $g=2$, the upper bound is due to Erd\H{o}s \& Tur\`{a}n~\cite{1941.Erdos.Turan} and the upper bound is due to
Singer~\cite{1938.Singer}. For $g=3$, the upper bound is Ruzsa's~\cite{1993.Ruzsa}, and the lower bound is from
$R(g+1,n)\le R(g,n)$. For $g>3$ the upper bound is a combination of Green~\cite{2001.Green} (for small $g$) and
Martin~\& O'Bryant~\cite{Martin.OBryant} (for large $g$). The lower bound for $g=4$ is due to Habsieger \&
Plagne~\cite{2002.Habsieger.Plagne}, the $g=6$ and $g=8$ bounds are due to
\cite{2002.Cilleruelo.Ruzsa.Trujillo}, and for all other even $g$ by Martin \& O'Bryant
\cite{2004.Martin.OBryant}. The lower bounds for odd $g>3$ are just $\sigma_2(2g)\le \sigma_2(2g+1)$.

The proof of $\sigma_2(2) = 1$ is succinct and elegant; we present it momentarily. The upper bound was found
initially in~\cite{1941.Erdos.Turan} and simplified in~\cite{1969.2.Lindstrom}. The lower bound was found in
\cite{1938.Singer} and simplified to this form in~\cite{1993.Ruzsa}.

\begin{thm} The largest Sidon subset of $[n]$ has $\sim \sqrt{n}$ elements, i.e.,
$\sigma_2 = 1$.
\end{thm}

\begin{proof} Let $1\le a_1<a_2< \dots < a_r \le n$ be a Sidon sequence, and
consider the differences (the parameter $u$ will be set to $\floor{n^{1/4}}$):
    \begin{center}
    $a_2-a_1, \quad a_3-a_2,\quad \dots,\quad a_r-a_{r-1}$ \\
    $a_3-a_1, \quad a_4-a_2,\quad  \dots, \quad a_r-a_{r-2}$ \\
    $\vdots$ \\
    $a_{u+1}-a_1, \quad a_{u+2}-a_2,\quad  \dots,\quad a_r-a_{r-u}$
    \end{center}
The $k$-th row contains $r-k$ differences, and since $\{a_i\}$ is a Sidon set, these differences are distinct.
That's a total of $\sum_{i=1}^u (r-i) = ru-\frac12 u(u+1)$ differences. The sum of all these differences is at
least
    $$\sum_{i=1}^{ru-u(u+1)/2} i = \frac12\left(ru-\tfrac12 u(u+1)\right)\left(ru-\tfrac12 u(u+1)+1\right).$$
On the other hand, the sum of the differences in the $k$-th row telescopes to
    $$\sum_{i=r-k+1}^r a_i - \sum_{i=1}^k a_i < kn, $$
and so the sum of all the differences is less than $\sum_{k=1}^u kn=nu(u+1)/2$. Comparing the upper and lower
bounds yields an inequality in terms of $n$, $r$, and $u=\floor{n^{1/4}}$. Calculus implies that
$r<n^{1/2}+n^{1/4}+1$, which in turn implies that $\sigma_2 \le 1$.

We now give Ruzsa's construction of a Sidon set contained in $p(p-1)$ with $p-1$ elements ($p$ is any odd
prime), from which it follows (using the elementary fact that the ratio between consecutive primes goes to 1)
that $\sigma_2 \ge 1$. Let $\theta$ be a primitive root modulo $p$, and consider the set $\cA$ of integers
$a_t$ ($1\le t < p-1$) defined by
    $$1\le a_t < p^2-p \quad \text{ and } \quad
    a_t \equiv t \pmod{p-1} \quad\text{ and } \quad
    a_t \equiv \theta^t \pmod{p}.$$
Now suppose, by way of contradiction, that there are three pairs $(a_{r_m},a_{v_m})\in \cA \times \cA$
satisfying $a_{r_m}+a_{v_m}= k $ (for some fixed $k\in\Z$). Each pair gives rise to a factorization modulo $p$
of
    $$x^2 - k x + \theta^k\equiv (x-a_{r_m})(x-a_{v_m}) \pmod{p},$$
using the critical relation $a_{r_m}a_{v_m}\equiv \theta^{r_m+v_m}=\theta^k \pmod{p}$. Factorization modulo $p$
is unique, so it must be that two of the three pairs are congruent modulo $p$, say
    \begin{equation}\label{modp:eq}
    a_{r_1}\equiv a_{r_2} \pmod{p}.
    \end{equation}
In this case, $\theta^{r_1} \equiv a_{r_1}\equiv a_{r_2} \equiv \theta^{r_2} \pmod{p}$. Since $\theta$ has
multiplicative order $p-1$, this tells us that $r_1 \equiv r_2 \pmod{p-1}$. Since $a_{r_m} \equiv r_m
\pmod{p-1}$ by definition, we have
    \begin{equation}\label{modp-1:eq}
    a_{r_1} \equiv a_{r_2} \pmod{p-1}.
    \end{equation}
Equations \eqref{modp:eq} and \eqref{modp-1:eq}, together with $a_{r_1}+a_{v_1} = k = a_{r_2}+a_{v_2}$ imply
that the pairs $(a_{r_1},a_{v_1}), (a_{r_2},a_{v_2})$ are identical, and so there are {\em not} three such
pairs. Thus, $\cA$ is a Sidon set. In particular, $\cA=\Ruzsa{p}{\theta}{1}$.
\end{proof}

A Sidon sequence $a_1<a_2<\dots<a_k$ is called ``short'' if $a_k-a_1$ is as small as possible.
Figure~\ref{fig:SidonSets} contains (up to reflection and translation) all of the short Sidon sequences with
$k\le 10$, and two of the short Sidon sequences with $k=11$; I don't know if there are more. It is rumored that
Imre Ruzsa has computed that $\min\{a_{12}-a_1\} = 85$, $\min\{a_{13}-a_1\}= 106$, $\min\{a_{14}-a_1\}=127$. It
is likely that for each $k$, there is short Sidon set with $a_1=0$ and $a_2=1$, while it seems likely that each
triple $(a_1,a_2,a_3)$ is the beginning of only finitely many short Sidon sequences.

\begin{figure}\label{fig:SidonSets}
    \begin{center}
        \begin{tabular}{ccc}
        $k$  & $\min\{a_k-a_1\}$ &             Witness              \\
        \hline
        2  &        1        &             \{0,1\}              \\
        3  &        3        &            \{0,1,3\}             \\
        4  &        6        &           \{0,1,4,6\}            \\
        5  &       11        &          \{0,1,4,9,11\}          \\
        &                 &          \{0,2,7,8,11\}          \\
        6  &       17        &       \{0,1,4,10,12,17\}         \\
        &                 &       \{0,1,4,10,15,17\}         \\
        &                 &       \{0,1,8,11,13,17\}         \\
        &                 &       \{0,1,8,12,14,17\}         \\
        7  &       25        &      \{0,1,4,10,18,23,25\}       \\
        &                 &      \{0,1,7,11,20,23,25\}       \\
        &                 &      \{0,1,11,16,19,23,25\}      \\
        &                 &      \{0,2,3,10,16,21,25\}       \\
        &                 &      \{0,2,7,13,21,22,25\}       \\
        8  &       34        &     \{0,1,4,9,15,22,32,34\}      \\
        9  &       44        &   \{0,1,5,12,25,27,35,41,44\}    \\
        10 &       55        &  \{0,1,6,10,23,26,34,41,53,55\}  \\
        11 &       72        & \{0,1,4,13,28,33,47,54,64,70,72\} \\
        &                 &  \{0,1,9,19,24,31,52,56,58,69,72\}  \\
        \end{tabular}
    \caption{Shortest Sidon sequences (from~\cite{2004.Martin.OBryant})}
    \end{center}
\end{figure}

Figure~\ref{fig:R(g,n)} gives the values of $n$ for which $R(g,n)-R(g,n-1)=1$. This table was computed using the
easiest algorithm and a small amount of time. I encourage the interested reader (or her students!) to extend it.

The construction given in~\cite{2002.Cilleruelo.Ruzsa.Trujillo} is optimized by choosing $x$ so that
$R(g,x)/\sqrt{gx}$ is maximized. This appears to happen (for each $g$) with a fairly small value of $x$; a
formula would lay to rest further optimization efforts (such as those in
\cites{2002.Habsieger.Plagne,2004.Martin.OBryant}). Even given a perfect optimization, however, it is unlikely
that this construction is optimal in any sense.

\begin{figure}\label{fig:R(g,n)}
\small
\begin{center}
$g$\vskip4pt $k$\quad
\begin{tabular}{|c||r|r|r|r|r|r|r|r|r|r|}
\hline
    &          2 &          3 &         4 &         5 &         6 & 7 & 8 & 9 &    10     &    11     \\ \hline\hline
 3  &          4 &            &           &           &           & & & & &           \\ \hline
 4  &          7 &          5 &           &           &           & & & & &           \\ \hline
 5  &         12 &          8 &         6 &           &           & & & & &           \\ \hline
 6  &         18 &         13 &         8 &         7 &           & & & & &           \\ \hline
 7  &         26 &         19 &        11 &         9 &         8 & & & & &           \\ \hline
 8  &         35 &         25 &        14 &        12 &        10 & 9 & & & &           \\ \hline
 9  &         45 &         35 &        18 &        15 &        12 & 11 & 10 &           &           &           \\ \hline
 10 &         56 &         46 &        22 &        19 &        14 & 13 & 12 &        11 &           &           \\ \hline
 11 &         73 &         58 &        27 &        24 &        17 & 15 & 14 &        13 &    12     &           \\ \hline
 12 &  $\leq 92$ &  $\leq 72$ &        31 &        29 &        20 & 18 & 16 &        15 &    14     &        13 \\ \hline
 13 & $\leq 123$ & $\leq 101$ &        37 &        34 &        24 & 21 & 18 &        17 &    16     &        15 \\ \hline
 14 & $\leq 140$ & $\leq 128$ &        44 &        40 &        28 & 26 & 21 &        19 &    18     &        17 \\ \hline
 15 & $\leq 163$ &            & $\leq 52$ & $\leq 47$ &        32 & 29 & 24 &        22 &    20     &        19 \\ \hline
 16 & $\leq 195$ &            &           &           &        36 & 34 & 27 &        24 &    22     &        21 \\ \hline
 17 &            &            &           &           & $\leq 42$ & $\leq 38$ &        30 &        28 &    24     &        23 \\ \hline
 18 &            &            &           &           &           & & 34 & 32 &    27     &        25 \\ \hline
 19 &            &            &           &           &           & & $\leq 38$ & $\leq 36$ &    30     &        28 \\ \hline
 20 &            &            &           &           &           & & & & 33     &        31 \\ \hline
 21 &            &            &           &           &           & & & & $\leq 37$ &        35 \\ \hline
 21 &            &            &           &           &           & & & & & $\leq 38$ \\ \hline
\end{tabular}
\caption{$\min\{n\colon R_2(g,n) \geq k\}$ (from~\cite{2004.Martin.OBryant}, extended by John Trono)}
\end{center}\end{figure}

\subsection{$h>2$}

\begin{figure}\label{fig:sigmah}
\begin{picture}(436,280)
    \put(428,11){$h$}
    \put(24,272){$\sigma_h$}
    \put(0,0){\includegraphics{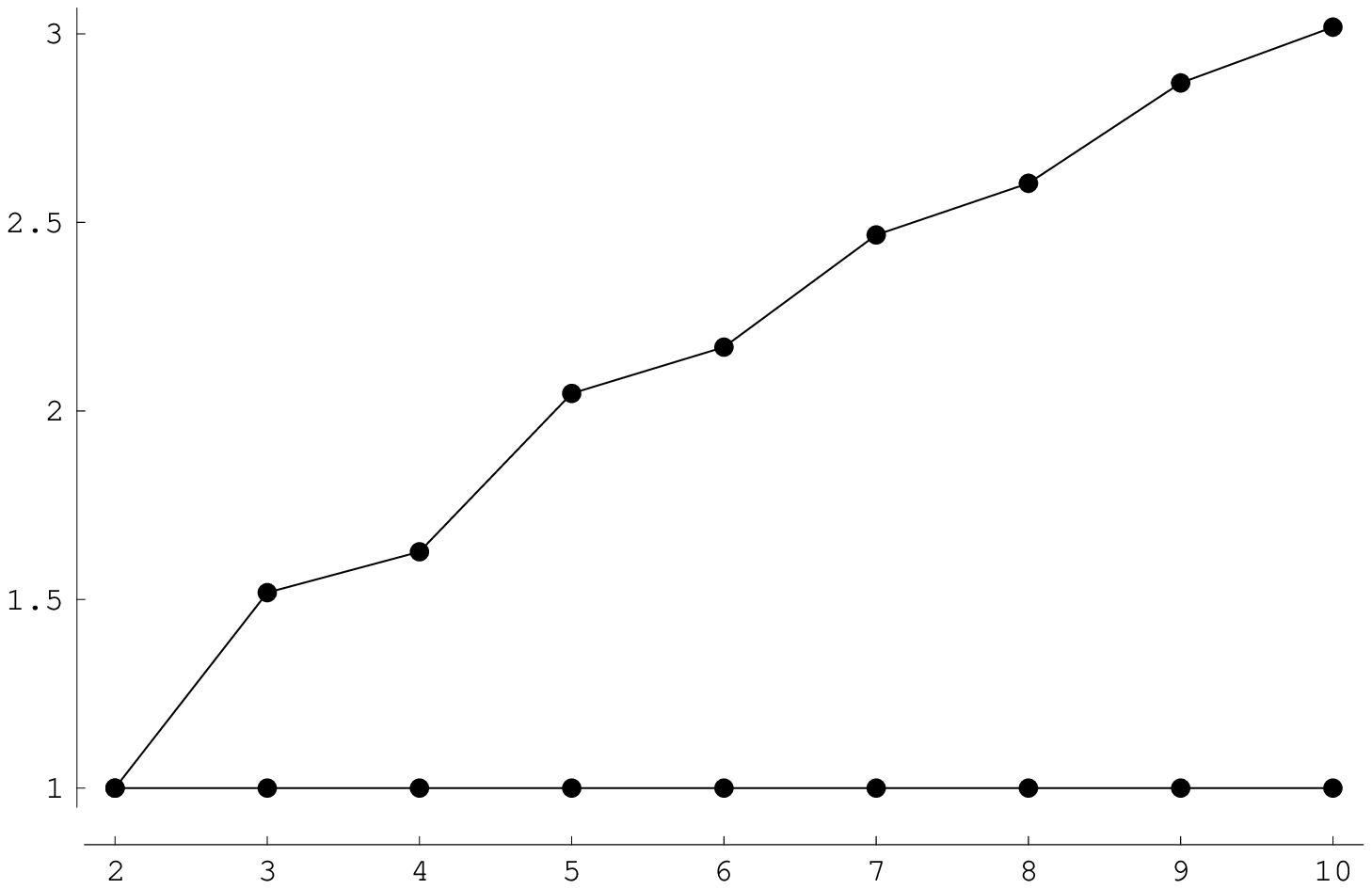}}
\end{picture}
\caption{The best known upper and lower bounds on $\sigma_h$.}
\end{figure}

Very little is understood about $B_h$ sets for $h>2$. The construction of Bose \& Chowla~\cite{1962.Bose.Chowla}
show that $\sigma_h\ge 1$, and we know that $\sigma_2=1$. Various improvements on the upper bound on $\sigma_h$
were given in~\cites{1961.Kruckeberg, 1969.3.Lindstrom, 1984.Dyachkov.Rykov, 1986.Shparlinski, 1993.Jia,
1994.Chen, 1996.1.Kolountzakis}

Cilleruelo~\cite{2001.Cilleruelo} proved that $\sigma_3 \le 1.576$, $\sigma_4 \le 1.673$, and for $h>75$,
    $$\sigma_h^h \le \frac 52 \left(\tfrac{15}{4}-\tfrac{5}{4\ceiling{h/2}}\right)^{1/4}
    \frac{(\ceiling{h/2}!)^2}{\sqrt{\ceiling{h/2}}} $$
for odd $h$, and
    $$\sigma_h^h \le \frac 52 \left(\tfrac{15}{4}-\tfrac{5}{2h}\right)^{1/4}
    ((h/2)!)^2 \sqrt{h/2}$$
for even $h$. He does give improvements for $h\in[5,74]$ also.

Ben Green~\cite{2001.Green} has improved these upper bounds and shown
    \begin{align*}
    \sigma_3 &\le (7/2)^{1/3} < 1.519 \\
    \sigma_4 &\le 7^{1/4} < 1.627  \\
    \sigma_h &\le \frac{1}{2e}\left( h + \frac32 \log h + o_h(\log h)\right).
    \end{align*}
Green's bound for $\sigma_h$ ultimately relies on the observation that if $X_i$ are independent random variables
taking values uniformly in a set of integers $\cA$, then the central limit theorem implies that
$X_1+X_2+\dots+X_h$ has a normal distribution (for large $h$). While Green's bound is not given in an effective
form, we presume that this could be done in a straightforward manner.

\subsection{Cyclic groups}

The earliest appearances of Sidon sets were as Monthly problems~\cites{1906.Veblen, 1906.Dickson} asking for
Sidon subsets of $\Z_n$. If $\cA$ is a $\Bhg{h}{g}(\bmod{ m})$ set and $\gcd(k,m)=1$, then so is $k \cA+r=\{ka+r
\colon a\in \cA\}$; we say that $\cA$ and $k\cA+r$ are equivalent. Veblen and Dickson asked for all inequivalent
Sidon subsets of $\Z_n$ for various $n$. The number of inequivalent modular Sidon sets that arise from the
construction of Bose \& Chowla is investigated in~\cites{1963.Halberstam.Laxton,1964.Halberstam.Laxton}.

The constructions of Ruzsa, Bose, and Singer all give {\em modular} Sidon sets. This would seem to be a more
symmetric setting, and so an easier setting. Unfortunately, while the constructions are naturally modular, the
upper bounds on $R_h(g,n)$ sets all seem to fundamentally rely on the {\em asymmetry} of $[n]$. Progress (beyond
the very little in~\cite{2004.Martin.OBryant}) on bounding $C_h(g,n)$ would be a significant contribution.
Here's what is known about $C_2(g,n)$:

\begin{thm} \label{C.upperbound}
    Let $q$ be a prime power, and let $k, g,f,x,y$ be positive integers with $k<q$.
    \renewcommand{\theenumi}{\roman{enumi}}
    \begin{enumerate}
    \item   $\binom{C_2(2,n)}{2} \le \floor{\frac n2}$, and in particular $C_2(2,n)\le \sqrt{n}+1$;
    \item   $C_2(3,n) \le \sqrt{n+9/2}+3$;
    \item   $C_2(4,n) \le \sqrt{3n} + 7/6$;
    \item   $C_2(g,n) \le \sqrt{gn}$ for even $g$;
    \item   $C_2(g,n) \le \sqrt{1-\tfrac1g}\sqrt{gn}+1$, for odd $g$.
    \item If $q$ is a prime, then $C_2(2k^2,q^2-q)) \geq k (q-1)$;
    \item $C_2(2k^2,q^2-1) \geq k q$;
    \item $C_2(2k^2,q^2+q+1) \geq kq+1$;
    \item If $\gcd(x,y)=1$, then $C_2(gf,xy)\geq C_2(g,x)C_2(f,y)$;
    \end{enumerate}
\end{thm}

\section{The size of infinite Sidon sequences}\label{sec.Size2}

Let $\cA\subseteq \Z$ be a Sidon sequence, and let $A(n)=\#(\cA \cap [n])$. St\"{o}hr~\cite{1955.Stohr} strengthens
an unpublished result of Erd\H{o}s and proved that
    $$\liminf_{n\to\infty} \frac{A(n)}{\sqrt{n/\log(n)}}\not=\infty$$
St\"{o}hr also gave Erd\H{o}s's proof that there is a Sidon sequence with
    $$\limsup_{n\to\infty} \frac{A(n)}{\sqrt{n}} \ge \frac 12;$$
this was improved by Kr\"{u}ckeberg~\cite{1961.Kruckeberg} to $1/\sqrt{2}$, and taken to $\Bhg{2}{g}$ sequences by
Cilleruelo \& Trujillo~\cite{2001.Cilleruelo.Trujillo}. In~\cite{1981.Ajtai.Komlos.Szemeredi} a Sidon sequence
is constructed with
    $$A(n)> 10^{-3} (n \log n)^{1/3}$$
for sufficiently large $n$. Erd\H{o}s asked~\cite{1982.Erdos} if there is a sequence $a_1<a_2<\dots$ with $a_k =
\littleo{k^{3-\epsilon}}$ for any positive $\epsilon$; Ruzsa constructed a Sidon sequence with $A(n)\sim
n^{\sqrt{2}-1+\littleo{1}}$.

Sheng Chen~\cite{1993.Chen} conjectures that if $\cA$ is a $B_{h}$ sequence with counting function $A(n)$, then
    $$\liminf_{n\to\infty} A(n) \left(\tfrac{\log n}{n}\right)^{1/h}$$
is finite. See also~\cites{1994.Helm, 1994.Jia, 1993.1.Helm, 1993.2.Helm, 1989.Nash, 1989.Jia, 1996.Chen,
1996.Helm, MR97d:11040}. Neither the results nor the conjectures have been extended to $\Bhg{h}{g}$ sequences.

\section{The distribution of $\cA$ and $\cA+\cA$}\label{sec.distribution}

If $\cA_n\subset [n]$ is a sequence of Sidon sets and $|\cA| \sim \sqrt{n}$, then $\cA_n$ becomes uniformly
distributed in $[n]$ (see~\cite{1991.Erdos.Freud}). Moreover, $\cA_n$ becomes uniformly distributed in
congruence classes modulo $m$ (for any fixed $m$) (see~\cites{1998.2.Lindstrom, 1999.Kolountzakis,
2002.Schoen}). It is shown in~\cite{1994.2.Erdos.Sarkozy.Sos} that if $\cA$ is any finite Sidon sequence, then
the sumset $\cA+\cA$ consists of at least $c|\cA|^2$ intervals (for an unspecified constant $c$). With care this
can be improved to $(|\cA|^2-|\cA|-1)/4$. Extensive computations indicate that if $|\cA| \sim \sqrt{n}$, then
$\cA+\cA$ consists of $\sim |\cA|^2/3$ intervals. Very little is known about the structure of $\cA+\cA$. Must
the size of the longest interval in $\cA+\cA$ go to infinity as $n$ does (with $|\cA|\sim\sqrt{n}$)?

Martin \& O'Bryant~\cite{Martin.OBryant} conjecture that $\Bhg{2}{g}$ sets with maximal size are uniformly
distributed. Current computations are insufficient to extend this conjecture to $\Bhg{h}{g}$ sets.

\section{Restricted Sidon sequences}\label{sec.SpecialSequences}

A well known conjecture~\cite{1994.Guy} states that the fifth powers $0, 1, 32, 243,\dots$ are a Sidon sequence.
Ruzsa~\cite{2001.Ruzsa} shows that there is an $\alpha$ such that $\{n^5+\floor{\alpha n^4}\colon n\ge n_0\}$ is
a Sidon set. Cilleruelo~\cite{1995.Cilleruelo} considered Sidon sequences all of whose terms are squares. Abbot
\cite{1990.Abbott} studied the size of Sidon sequences contained in an arbitrary set of integers with
cardinality $n$ (there is one with size at least $\frac 2{25} \sqrt{n}$).

Lindstr\"{o}m~\cites{1969.1.Lindstrom,1972.Lindstrom} initiated the study of $\Bhg{h}{g}$ sets in group $\Z^d$.

Erd\H{o}s~\cite{1980.Erdos} investigated Sidon subsets of $\Bhg{2}{g}$ sets. Alon \& Erd\H{o}s
\cite{1985.Alon.Erdos} considered the problem of decomposing a finite $\Bhg{2}{g}$ into $B_2$ sets; how many
$B_2$ sets are needed?

\section{Generalizations}

The graph theorist's analogs of Sidon sets are magic and harmonious labelings. See
\cites{1980.1.Graham.Sloane,1980.2.Graham.Sloane}.

The bigger setting for the Sidon question is that of sets which avoid a linear form. Let $L$ be an $m\times n$
matrix of integers. The set $\cA$ is said to avoid $L$ if there are no nontrivial solutions to
$L\vec{a}=\vec{0}$ with $\vec{a}=(a_1,a_2, \dots ,a_n)^T, a_k\in \cA$. Sidon sets are the special case
$L=[1,1,-1,-1]$, sets with 3-term arithmetic progressions are the case $L=[1,-2,1]$, etc. It is surprising that
in such a general setting significant and strong results can be found. Precisely this is done
in~\cites{1972.Komlos.Sulyok.Szemeredi, 1975.Komlos.Sulyok.Szemeredi,1993.Ruzsa}.

\section{Other open questions}\label{sec.OpenQuestions}

If $\cA^\ast$ is bounded, is it necessarily 0 infinitely often? This is equivalent to a USD 500 question of
Erd\H{o}s~\cite{1989.Erdos}: If every positive integer can be written as a sum of two elements of
$\cB\subseteq\N$ (i.e., $\cB$ is a base), must the number of ways to do so go to infinity?

Is there an anti-Freiman theorem: If $|\cA+\cA|/|\cA|^2 > c$, then $\cA$ has a large $\Bhg{2}{g}$ subset, where
$g$ and `large' depend only $c$? This question may be related to the quasi-Sidon sets discussed
in~\cites{1991.Erdos.Freud, Pikhurko}.

\end{document}